\documentclass[letterpaper, 10 pt, conference]{ieeeconf}  

\IEEEoverridecommandlockouts                              
\usepackage{graphics} 
\usepackage{epsfig} 
\usepackage{amsmath} 
\usepackage{amssymb}  
\usepackage{caption}
\usepackage{subcaption}

\title{\LARGE \bf
Synchronization of Heterogeneous Kuramoto Oscillators with Arbitrary Topology
}

\author{Authors}
\author{Andrey Gushchin, Enrique Mallada, and  Ao Tang 
\thanks{Andrey Gushchin is with the Center for Applied Mathematics,
        Cornell University, Ithaca, NY 14850, USA
        {\tt\small avg36@cornell.edu}}
\thanks{Enrique Mallada is with the Center for the Mathematics of Information, Department of Computing + Mathematical Sciences, Caltech, Pasadena,
        CA 91125, USA
        {\tt\small mallada@caltech.edu}}%
\thanks{Ao Tang is with the School of Electrical and Computer Engineering, Cornell University,
        Ithaca, NY 14850, USA
         {\tt\small atang@ece.cornell.edu}}%
}

\begin{document}

\maketitle
\thispagestyle{empty}
\pagestyle{empty}

\begin{abstract}
We study synchronization of coupled Kuramoto oscillators with heterogeneous inherent frequencies and general underlying connectivity. We provide conditions on the coupling strength and the initial phases which guarantee the existence of a Positively Invariant Set (PIS) and lead to synchronization. Unlike previous works that focus only on analytical bounds, here we introduce an optimization approach to provide a computational-analytical bound that can further exploit the particular features of each individual system such as topology and frequency distribution. Examples are provided to illustrate our results as well as the improvement over previous existing bounds.
\end{abstract}

\section{INTRODUCTION}


The study of synchronization of coupled oscillators has attracted the attention of widely diverse research disciplines such as neuroscience~\cite{c10,c9,c11}, physics~\cite{c13,c12}, mathematics~\cite{c5} and engineering~\cite{c24,c26,c27}. 
Since the seminal works of Winfree~\cite{c9} and Kuramoto~\cite{c6}, the Kuramoto model has served as a canonical model for synchronization that can capture a quite rich dynamic behavior including multiple equilibria, limit cycles, and even chaos.

There are two main properties that characterize its behavior.
The first one is the coupling function, which is a trigonometric $\sin()$ function in the case of the Kuramoto model. However, a broader class of the coupling functions has also been studied~\cite{c23,c7,c22}. 
The second property, and perhaps the most important one, is the interconnection topology. The most popular assumption is that all oscillators are connected to each other, which corresponds to a fully connected graph~\cite{c1,c2}. Although a much more general approach is to study the systems of oscillators with arbitrary underlying topology~\cite{c20,c3,c5,c19}. 

Due to its complex behavior, several assumptions are usually made to make the study tractable. For example, one can make the number of oscillators go to infinity, and use statistical mechanics tools to characterize its convergence. Or one can assume that all oscillators have equal intrinsic frequencies, and therefore form a (gradient) system of homogeneous oscillators~\cite{c7} that has globally convergent properties.
Alternatively, as we do in this article, one may let the frequencies take distinct values~\cite{c15,c1,c2,c18,c17,c16} and characterize sufficient conditions for synchronization.

In this paper, we consider a system of finite number of heterogeneous Kuramoto oscillators with general underlying topology. 
We show that when certain conditions on the coupling strength and initial phases are satisfied, the trajectories are bounded, which for Kuramoto oscillators also implies synchronization. The most relevant previous work is~\cite{c5}, where the authors studied the same setup. Here, we build upon~\cite{c5} to obtain less restrictive conditions for synchronization. 
In particular, by using a novel computational-analytical approach, we are able to further exploit the particular features of each problem instance and outperform existing results.
Several examples are used to illustrate our findings and characterize the scaling behavior of our conditions.


\section{PROBLEM FORMULATION}
We consider a system of oscillators which are described by a Kuramoto model, i.e. the behavior of each oscillator is governed by the following equation:
\begin{equation}\label{eq:1}
\dot \phi_{i}=w_{i}+\frac{K}{n}\sum_{k \in N_{i}} \sin(\phi_{k}-\phi_{i}),
\end{equation}
where $N_{i}$ is a set of oscillators connected to oscillator $i$, i.e. the set of its neighbors, $K$ is the coupling strength, which assumed to be the same for all connections, and $n$ is the total number of oscillators in the system. We also assume that the oscillators are heterogeneous in intrinsic frequencies. This means that their intrinsic frequencies $w_{i}$ are not necessary equal. Frequencies, however, do not change their values with time, so each $w_{i}$ is  a constant.

 In this article we study frequency synchronization of the system \eqref{eq:1}. The oscillators achieve synchronization if $\dot \phi_{1}(t) = \dot \phi_{2}(t) = \dots = \dot \phi_{n}(t)=\dot \phi$ as $t \rightarrow \infty $, where $\dot \phi$ is a constant common phase velocity.

It is easy to show that this common phase velocity $\dot \phi$ is an average sum of intrinsic frequencies of the oscillators. That is,~
$$  \dot \phi = \frac{\sum\limits_{k=1}^{n}w_{i} }{n}.    $$

Indeed, when $\dot \phi_{1} = \dot \phi_{2} = \dots = \dot \phi_{n}$,  adding up all the equations of \eqref{eq:1} gives: $(\dot \phi_{1} + \dot \phi_{2} +\dots +\dot \phi_{n} ) = w_{1}+ w_{2}+ \dots+~w_{n}$, because each $\frac{K}{n} \sin(\phi_{k}-\phi_{i})$ is added to $\frac{K}{n}\sin(\phi_{i}-\phi_{k})$ and gives zero.

We denote the average natural frequency by $\bar{w}\triangleq\frac{\sum\limits_{k=1}^{n}w_{i}}{n}$, and define the deviations of the natural frequencies by\\ $\bar{w}_{i}\triangleq w_{i}-\bar{w}$, where $i=1,\dots,n$.
From now on we will study the following system instead of system \eqref{eq:1}:

\begin{equation}\label{eq:2}
\dot \phi_{i}=\bar{w}_{i}+\frac{K}{n}\sum_{k \in N_{i}} \sin(\phi_{k}-\phi_{i}).
\end{equation}

Each limit cycle of system \eqref{eq:1} is an equilibrium of \eqref{eq:2}. Therefore, we will focus on finding conditions when system \eqref{eq:2} converges to an equilibrium, i.e. when $\dot{\phi_{i}}=0 \; \forall i=1,\dots,n$. Due to the rotational invariance of \eqref{eq:2}, we can shift without loss of generality, all the initial phases $\phi_{1}^{0}, \dots, \phi_{n}^{0}$ by the same value $\sum\limits_{i=1}^{n}\phi_{i}^{0}/n$ so that their sum becomes equal to zero: $\sum\limits_{i=1}^{n}\phi_{i}^{0}=0$. 
Furthermore, since the phase average remains the same (for system \eqref{eq:2}: $\dot{\phi_{1}}+\dots+\dot{\phi_{n}}=0$), condition $\sum\limits_{i=1}^{n}\phi_{i}^{t}=0$ will be satisfied $\forall t\geq 0$, where $\vec{\phi}^{t}$ are the trajectories of system \eqref{eq:2}. In the rest of this article we will assume that the phase values sum up to zero at each time $t \geq 0$.

In this article we show synchronization of oscillators by providing a Lyapunov function and using LaSalle's Invariance Theorem \cite{c21}. When all the intrinsic frequencies are equal, i.e. deviations $\bar{w}_{1}= \dots= \bar{w}_{n}=0$, the oscillators are called homogeneous, and the following Lyapunov function can be employed:
$$ V_{0}(\vec \phi)= -\frac{K}{n}\sum_{ij \in E, i<j} \cos(\phi_{i}-\phi_{j}),$$
where $\vec\phi\in\mathbb{R}^{n}$ and $E$ is the edge set of a given graph. It is easy to check that
$$ \dot{V_{0}}(\vec \phi) =-\sum\limits_{i=1}^{n}\dot{\phi}_{i}^{2} \leq 0.$$
%
%
Thus, since the function $V_{0}(\vec \phi)$ is $2\pi$-periodic on each element $\phi_i$ in $\mathbb R$, it is also well-defined on a n-dimensional torus ($\mathbb T^n$), which is compact. Thus, applying the LaSalle's Invariance Theorem (on $\mathbb T^n$) guarantees synchronization of the oscillators.\footnote{See the proof of Proposition 1 for an example on how to apply LaSalle's Invariance Theorem.}

When the intrinsic frequencies are not equal, we have a system of heterogeneous oscillators, and we still can provide a potential function for this case:
 $$ V(\vec \phi)\triangleq -\sum_{k=1}^{n} (\bar{w}_{k}\phi_{k})-\frac{K}{n}\sum_{ij \in E, i<j} \cos(\phi_{i}-\phi_{j}). $$
We can check again that the time derivative of this function is also non-positive and equal to zero only at an equilibrium, i.e. when the frequencies are synchronized.

The problem here is that function $V(\vec \phi)$ is not bounded from below and cannot be defined on the $n$-dimensional torus $\mathbb{T}^n$, and therefore, we are not able to apply directly the LaSalle's Invariance Theorem. However, if we show that the trajectories $\vec \phi^t\in\mathbb R^n$ of \eqref{eq:2} are bounded, then the function $V(\vec \phi)$ is bounded as well, and hence synchronization follows.

Therefore, the main goal of this study is to find the conditions that guarantee that trajectories are bounded.
We will achieve this by finding a compact Positively Invariant Set (PIS) for the oscillators' phases. That is, a compact (closed and bounded) set such that if the system's initial conditions are within this set, the trajectories will remain in the set. 

The next section shows that when some conditions are met, such PIS exists, and therefore system \eqref{eq:2} will converge to the set of equilibria.

%

\section{MAIN RESULTS}
This section is organized as follows. We first formulate in Proposition 1 a general sufficient condition for boundedness of the trajectories that leads to synchronization of system \eqref{eq:2}. We then provide two solutions that guarantee fulfillment of Proposition 1. Our first solution, described in subsection ${\it{B}}$, contains explicit requirements on the coupling strength and initial oscillators' phases. This solution is further refined using computational tools in subsection ${\it{C}}$.

\subsection{Preliminary Results}

We will denote maximum and minimum phase values at time $t$ by $\phi_{max}^{t} \triangleq \max\limits_{i}\phi_{i}^{t}$ and $\phi_{min}^{t} \triangleq \min\limits_{i}\phi_{i}^{t}$, where $\phi_{i}^{t}$ is a phase of oscillator $i$ at time $t$.
Let $D_{t}$ be defined as a maximum phase difference between two oscillators at time $t$ $(t\geq 0)$, i.e. 
$$D_{t}~\triangleq~ \phi_{max}^{t}-\phi_{min}^{t},$$
then $\phi_{min}^{t} \leq \phi_{i}^{t} \leq \phi_{max}^{t}$ ($\forall i = 1, \dots, n$). In other words, each phase lies between the minimum and maximum phases $\phi_{min}^{t}$ and $\phi_{max}^{t}$.
The maximum initial (at time $t=0$) pairwise phase difference is denoted by $D_{0}$:  
$$D_{0}~=~ \phi_{max}^{0}-\phi_{min}^{0}.$$
If we can show that the maximum phase difference is always bounded, i.e. if $D_{t} \leq D \; \forall \;t \geq 0$, where $D$ is a constant satisfying $D_{0} \leq D <\infty$, then the trajectories will be also bounded since the phase average remains the same. The PIS therefore is defined through the maximum phase difference that is bounded by the value of $D$, i.e.
$$\text{PIS}\triangleq\{\vec{\phi}\in \mathbb{R}^{n}: \; \max_{i,j}|\phi_i-\phi_j|\leq D, \; \sum\limits_{i=1}^{n}\phi_{i}=0\},$$
which is obviously a compact.

We can now formulate a general condition that is sufficient to guarantee that the maximum phase difference is always bounded by a constant $D$ and thus the trajectories are also bounded.

{\bf{Proposition 1}} {\it{ If $D$ is a constant satisfying $D_{0} \leq D <\infty$, and for all times $t\geq 0$ such that $D_{t}=\phi_{max}^{t}-\phi_{min}^{t} = D$, the following condition is satisfied:
\begin{equation}\label{eq:3}
\begin{split}
&\dot{\phi}_{k}^{t}-\dot{\phi}_{l}^{t} = \bar{w}_{k}-\bar{w}_{l} \\ 
&- \frac{K}{n}\sum\limits_{i\in N_{k}}\sin(\phi_{k}^{t} - \phi_{i}^t) -\frac{K}{n}\sum\limits_{j\in N_{l}}\sin(\phi_{j}^{t} - \phi_{l}^t) \leq 0,  \\ 
\end{split} 
\end{equation}
for every two oscillators $k$ and $l$ such that $\phi_{k}^{t}=\phi_{max}^{t}$ and $\phi_{l}^{t}=\phi_{min}^{t}$, then the maximum phase difference is bounded by $D$, i.e. $D_{t}\leq D$ for all $t\geq 0$, trajectories of system \eqref{eq:2} are bounded, and system \eqref{eq:2} achieves frequency synchronization.}}

\begin{proof}
Condition \eqref{eq:3} says that when the maximum phase difference achieves value $D$, it cannot grow anymore and thus does not exceed $D$. That's why the maximum phase difference will be always bounded by $D$ if \eqref{eq:3} is satisfied, and because the phase average remains the same, it also ensures that the trajectories of system \eqref{eq:2} are bounded in $\mathbb{R}^{n}$. Since function $V(\vec{\phi})$ is well-defined in $\mathbb{R}^{n}$, we can now apply LaSalle's Invariance Theorem to guarantee that each solution of \eqref{eq:2} approaches the nonempty set $\{\dot V\equiv 0\}=\{\dot{\phi}_{i}=0, 1\leq~i\leq~n\}$, and system \eqref{eq:2} achieves frequency synchronization.
\end{proof}
It is possible that when $\phi_{max}^{t}-\phi_{min}^{t} = D$, several oscillators have phase values equal to $\phi_{max}^{t}$ or $\phi_{min}^{t}$. In this case condition \eqref{eq:3} should be satisfied for each pair of oscillators with a phase difference equal to $D$. 

Condition \eqref{eq:3} is very general by itself and difficult to check.
In the rest of this section we derive two conditions -- analytic and optimization-based -- that guarantee condition \eqref{eq:3}. These conditions contain requirements on the coupling strength and initial phases of oscillators that can be verified for each given system. 
We now introduce some additional notation. 

Let
$$\mathcal{E}_{t}(\phi) \triangleq \sum\limits_{i=1}^{n}\bigl(\phi_{i}^{t}\bigl)^{2},$$
i.e. $\mathcal{E}_{t}(\phi)$ is the squared Euclidean norm of a vector of phases at time $t$. For simplicity we will use symbol $\mathcal{E}_{t}$ instead of $\mathcal{E}_{t}(\phi)$. At initial time $t=0$ the value of this function is denoted by $\mathcal{E}_{0}$.

The Euclidean norm of a vector of the natural frequencies deviations is defined as $\sigma(\bar{w})$:
$$\sigma(\bar{w}) \triangleq \sqrt{\sum\limits_{i=1}^{n}(\bar{w}_{i})^2}.$$

Let the topology of a given system be defined by an undirected graph $G=(V,E)$ with a set of nodes $V$ such that $|V|=n$, and with an edge set $E$. By $E^{c}$ we denote the set $E_{comp} \setminus E:$
$$E^{c}=E_{comp} \setminus E,$$
where $E_{comp}$ is the set of $\frac{n(n-1)}{2}$ edges of a complete graph with $n$ nodes. We use $\delta$ to denote the minimum nodal degree of a graph.

The following two lemmas are based on the results from \cite{c5} and will be used in the next two subsections; their proofs are provided in appendix. 


{\bf{Lemma 1}} {\it{If $\sum\limits_{i=1}^{n}\phi_{i}^{t}=0$, then

 $$ L \cdot n \cdot \mathcal{E}_{t} \leq \sum\limits_{(i,j) \in E}|\phi_{i}^{t}-\phi_{j}^{t}|^{2} \leq n\cdot \mathcal{E}_{t}, $$
where  $L \triangleq \frac{1}{1+\sum\limits_{(k,l) \in E^{c}}dist(k,l)}$.}}

Distance $dist(k,l)$ between two nodes $k, l \in V$  in a graph $G=(V,E)$ is the number of edges in a shortest path between these two nodes.



{\bf{Lemma 2}} {\it{ If $D_{t} \leq D<\pi \; \forall t \in [0, T]$, then function $\mathcal{E}$ satisfies the following differential inequality on $[0,T]$:
\begin{equation}\label{eq:4}
\frac{d}{dt}\mathcal{E} \leq 2\sigma(\bar{w})\cdot  \sqrt{\mathcal{E}}-2K\cdot L\cdot \Bigl(\frac{\sin D}{D}\Bigl)\cdot \mathcal{E}, 
\end{equation}
and if in addition
\begin{equation}\label{eq:5}
 K\geq \frac{\sigma(\bar{w})\cdot D}{\sqrt{\mathcal{E}_{0}}\cdot L\cdot \sin D} \;, 
 \end{equation}
then $\mathcal{E}_{t}$ will be upper bounded by $\mathcal{E}_{0}$: $\mathcal{E}_{t} \leq \mathcal{E}_{0} \; \forall t\in [0, T]$.}}
It is assumed that not all initial phases are equal to zero, so that $\mathcal{E}_{0}>0$ in condition \eqref{eq:5}. 

Lemma 2 states that if the trajectories of system \eqref{eq:2} stay in a PIS defined by $D_{t}\leq D<\pi$, and the coupling strength $K$ satisfies condition \eqref{eq:5}, then function $\mathcal{E}$ does not exceed its initial value.

\subsection{Analytic Synchronization Condition}
The main result of this subsection is Theorem 1 which contains requirements on the initial phases and coupling strength such that all conditions of Proposition 1 are satisfied and thus system \eqref{eq:2} achieves frequency synchronization.
 \begin{figure}[t]
      \centering
      \includegraphics[scale=0.5]{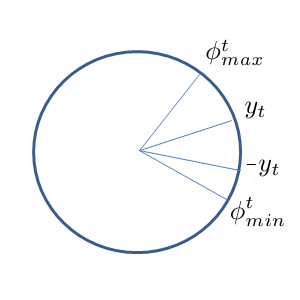}
      \caption{$\phi_{max}^{t}, \phi_{min}^{t}$, $y_{t}$ and $-y_{t}$ at time $t$ when $D_{t}= D.$}
      \label{five}
   \end{figure}

$\bf{Theorem \; 1}$ {\it {If $D$ is a constant satisfying $0<D_{0} \leq D \leq \frac{\pi}{2}$; $\mathcal{E}_{0}<\frac{3}{4}D^{2}$; $K$ satisfies \eqref{eq:5} and
\begin{equation}\label{eq:6}
K \geq \frac{n\cdot |\bar{w}_{i}-\bar{w}_{j}|}{2\delta \cdot \sin\Bigl( \frac{D}{2}-\sqrt{\mathcal{E}_{0}-\frac{D^{2}}{2}}\Bigl)}  \; \forall i,j,
\end{equation}
where $1 \leq i,j \leq n$, then there exists a PIS defined by $D_{t}\leq D$ for all $t \geq 0$, and system \eqref{eq:2} achieves frequency synchronization. When $\mathcal{E}_{0}<\frac{D^{2}}{2}$, only condition \eqref{eq:5} for $K$ is required. }}

\begin{proof}
%
We first consider the case when $\mathcal{E}_{0}\geq\frac{D^{2}}{2}$.
Assume that at time moment $t\geq 0$: $D_{t}=D$, and that before this moment $t$ the maximum phase difference has never exceeded $D$. Then if \eqref{eq:5} holds, function $\mathcal{E}_{t}$ has not exceeded its initial value $\mathcal{E}_{0}$. 

We define:
 $$y_{t} \triangleq \sqrt{\mathcal{E}_{t} - (\phi_{min}^{t})^{2}-(\phi_{max}^{t})^{2}}.$$
Suppose that $k$ and $l$ are two oscillators with $\phi_{k}^{t}=\phi_{max}^{t}$ and $\phi_{l}^{t}=\phi_{min}^{t}$.
By the definition of $y_{t}$:  
\begin{equation}\label{eq:7}
-y_{t} \leq \phi_{i}^{t} \leq y_{t} 
\end{equation}
$\forall i=1,\dots,n,$ such that $i \neq k$, $i \neq l$. On Fig. 1 we plotted phases $\phi_{min}^{t}$, $\phi_{max}^{t}$ at time $t$ with $D_{t}=\phi_{max}^{t}-\phi_{min}^{t} = D$, and values of $\pm y_{t}$.

If at time $t$ $D_{t}=D$, then because $\sum\limits_{i=1}^{n}\phi_{i}^{t}=0$, one of the following possibilities takes place: $\phi_{max}^{t} > D/2$, or $\phi_{min}^{t} < -D/2$, or $\phi_{max}^{t}=D/2, \;\phi_{min}^{t}=-D/2$. 
Let  
\begin{equation}\label{eq:phi-max/min}
\phi^{t}_{max} = D/2+d_t \quad\text{ and }\quad \phi^{t}_{min} = -D/2+d_t,
\end{equation}
then
\begin{equation*}
\begin{split}
\mathcal{E}_{t} &= (D/2+d_t)^{2} + (-D/2+d_t)^{2}+y_{t}^{2} \\
&= D^{2}/2 + 2d_t^{2}+y_{t}^{2} \leq \mathcal{E}_{0}.
\end{split}
\end{equation*}
Therefore,
\begin{equation}\label{eq:y_t-bound}
y_{t}^{2} \leq \mathcal{E}_{0}-\frac{D^{2}}{2}-2d_t^{2},
\end{equation}
and
\begin{equation}\label{eq:d_t-bound}
 d_t^{2}  \leq \frac{\mathcal{E}_{0}-\frac{D^{2}}{2}}{2}.
\end{equation}

We want to show that at time $t$ the maximum phase difference does not start to increase by showing that condition \eqref{eq:3} of Proposition 1 is satisfied, i.e.
$\dot{\phi}_{k}^{t}-\dot{\phi}_{l}^{t} \leq 0$ for every two oscillators $k$ and $l$ with $\phi_{k}^{t}=\phi_{max}^{t}$ and $\phi_{l}^{t}=\phi_{min}^{t}$. 

Using equation \eqref{eq:7}, \eqref{eq:phi-max/min} and the fact that $D\leq \frac{\pi}{2}$ we have 
\begin{align}
\frac{\pi}{2}&\geq\phi_k^t-\phi_i^t \geq  \phi_k^t-y_t = \frac{D}{2}+d_t-y_t\label{eq:inequality-0}\\
&\geq \frac{D}{2}-\frac{D}{2}- \sqrt{ \mathcal{E}_{0}-\frac{D^{2}}{2}-2d_t^{2}} \label{eq:inequality-1}\\
&> -\sqrt{\frac{3D^2}{4}-\frac{D^{2}}{2}}\label{eq:inequality-2}=-\frac{D}{2}\geq-\frac{\pi}{4}
\end{align}
where \eqref{eq:inequality-1} follows from \eqref{eq:y_t-bound} and $d_t\geq -\frac{D}{2}$, and \eqref{eq:inequality-2} from $-\sqrt{\cdot}$ being decreasing and the theorem assumption $\mathcal{E}_{0}<\frac{3}{4}D^{2}$.
Similarly, we have 
\begin{equation}\label{eq:inequality-3}
\frac{\pi}{2}\geq \phi_j^t-\phi_l^t\geq -y_t -\phi_l^t\geq-\frac{\pi}{4}.
\end{equation}

Now, from condition \eqref{eq:3} of Proposition 1 we have
\begin{align}
&\dot{\phi}_{k}^{t}-\dot{\phi}_{l}^{t} = \bar{w}_{k}-\bar{w}_{l} \nonumber\\
&- \frac{K}{n}\sum\limits_{i\in N_{k}}\sin(\phi_{k}^{t} - \phi_{i}^t)-\frac{K}{n}\sum\limits_{j\in N_{l}}\sin(\phi_{j}^{t} - \phi_{l}^t) \nonumber\\
&\leq \bar{w}_{k}-\bar{w}_{l}-\frac{\delta \cdot K}{n}\Bigl(\sin(\phi_{k}^{t}-y_{t})+\sin(-\phi_{l}^{t}-y_{t} )\Bigl)\label{eq:inequality}\\
&= \bar{w}_{k}-\bar{w}_{l}-\frac{\delta \cdot K}{n}\Bigl(\sin(\phi_{max}^{t}-y_{t})+\sin(-\phi_{min}^{t}-y_{t} )\Bigl).\nonumber
\end{align}
where inequality \eqref{eq:inequality} follows from \eqref{eq:inequality-0}-\eqref{eq:inequality-2}, \eqref{eq:inequality-3}, the fact that $\sin()$ is an increasing function on $(-\frac{\pi}{2}, \frac{\pi}{2})$, and because $0\leq (\phi_{k}^{t}-\phi_{i}^{t})\leq \frac{\pi}{2}$, $0\leq (\phi_{j}^{t}-\phi_{l}^{t})\leq \frac{\pi}{2} \; \forall i$ and $\forall j$.
%

We will now show that
\begin{equation}\label{eq:8}
\begin{split}
\sin(\phi^{t}_{max}-y_{t})&+\sin(-\phi^{t}_{min}-y_{t})\\
&\geq 2\sin\Bigl(\frac{D}{2}-\sqrt{\mathcal{E}_{0}-\frac{D^{2}}{2}}\Bigl) >0
\end{split}
\end{equation}
where the last inequality holds because $\sqrt{\mathcal{E}_{0}-\frac{D^{2}}{2}}<\frac{D}{2}.$

For simplicity we will introduce the following notation:
\begin{equation}
\label{alpha}
\alpha_{t} \triangleq \sqrt{\mathcal{E}_{0}-\frac{D^{2}}{2}-2d_t^{2}}.
\end{equation}
Since due to \eqref{eq:y_t-bound},  $y_{t} \leq \alpha_{t}$, from \eqref{eq:phi-max/min} and  \eqref{eq:inequality-0}-\eqref{eq:inequality-2}, \eqref{eq:inequality-3} we have:
\begin{equation}\label{eq:18}
\frac{\pi}{2}\geq \phi^{t}_{max}-y_{t} \geq \frac{D}{2}+d_t-\alpha_{t}\geq -\frac{\pi}{4}, 
\end{equation}
 and
\begin{equation}\label{eq:19}
\frac{\pi}{2}\geq -y_{t}-\phi^{t}_{min} \geq \frac{D}{2}-d_t-\alpha_{t}\geq -\frac{\pi}{4}. 
\end{equation}
Consider a function $f(d)$ of one scalar argument $d$:
\begin{equation*}
\begin{split}
f(d) = \sin\Bigl( \frac{D}{2}+d-\alpha\Bigl)+ \sin\Bigl( \frac{D}{2}-d-\alpha \Bigl),
\end{split}
\end{equation*}
where $d$ satisfies \eqref{eq:d_t-bound} and $\alpha$ is a function of $d$ defined by \eqref{alpha}. We will now show that $f(d)\geq 2\sin\Bigl(\frac{D}{2}-\sqrt{\mathcal{E}_{0}-\frac{D^{2}}{2}}\Bigl)$.
When $d=0$ we get an equality:
\begin{equation*}
\begin{split}
f(0)=2\sin\Bigl( \frac{D}{2}-\alpha\Bigl)
=2\sin\Bigl(\frac{D}{2}-\sqrt{\mathcal{E}_{0}-\frac{D^{2}}{2}}\Bigl).
\end{split}
\end{equation*}
Assume first that $d\geq 0$. Now it is enough to show that the derivative of this function is positive: $f^{'}(d)>0$ for all $0 \leq d  \leq \sqrt{\frac{\mathcal{E}_{0}-\frac{D^{2}}{2}}{2}}$.
It can be verified, that
\begin{equation*}
\begin{split}
f'(d) = &-2\sin\Bigl(\frac{D}{2}-\alpha\Bigl)\cdot \sin d \\
&+\frac{4d}{\alpha}\cdot \cos d \cdot \cos\Bigl(\frac{D}{2}-\alpha\Bigl).
\end{split}
\end{equation*}
Notice, that $\frac{D}{2}-\alpha \in (0, \frac{D}{2}] \in (0,\frac{\pi}{4}]$, and 
$\sin\Bigl(\frac{D}{2}-\alpha\Bigl)>0$, $\cos\Bigl(\frac{D}{2}-\alpha\Bigl)>0$ and $0<\tan\Bigl(\frac{D}{2}-\alpha\Bigl)\leq 1.$

 Since $0 \leq d  \leq \sqrt{\frac{\mathcal{E}_{0}-\frac{D^{2}}{2}}{2}}$, $\sin d\geq 0, \cos d > 0$ and $\sin d \leq d$.
Therefore,
\begin{equation*}
\begin{split}
f'(d) &\geq -2d\cdot \sin\Bigl(\frac{D}{2}-\alpha\Bigl)+\frac{4d}{\alpha}\cdot \cos d\cdot \cos\Bigl(\frac{D}{2}-\alpha\Bigl) \\
&= 2d\cdot \cos\Bigl(\frac{D}{2}-\alpha\Bigl)\cdot \Bigl( -\tan\Bigl(\frac{D}{2}-\alpha\Bigl)+\frac{2\cos d}{\alpha}\Bigl).
\end{split}
\end{equation*}
Now, since $-1 \leq -\tan\Bigl(\frac{D}{2}-\alpha\Bigl) < 0$, it is sufficient to show that $2\cos d > \alpha$. Indeed, since $\mathcal E_0<\frac{3D^2}{4}$, $\cos d \geq \cos\Bigl( \sqrt{\frac{\mathcal{E}_{0}-\frac{D^{2}}{2}}{2}}\Bigl)>\cos(\frac{D}{2\sqrt{2}}) \geq \cos(\frac{\pi}{4\sqrt{2}})$, which means that $2\cos d> 2\cos(\frac{\pi}{4\sqrt{2}})>1.$ On the other hand, $\alpha \leq \sqrt{\mathcal{E}_{0}-\frac{D^{2}}{2}} \leq \frac{D}{2} \leq \frac{\pi}{4}<1.$

For $ -\sqrt{\frac{\mathcal{E}_{0}-\frac{D^{2}}{2}}{2}} \leq d \leq 0$ the proof is similar. Therefore, $f(d)\geq 2\sin\Bigl(\frac{D}{2}-\sqrt{\mathcal{E}_{0}-\frac{D^{2}}{2}}\Bigl)$, and \eqref{eq:8} holds because of \eqref{eq:18} and  \eqref{eq:19}.
Notice that even though it is possible that one of the $\sin$ functions $\sin(\phi^{t}_{max}-y_{t})$ or $\sin(-\phi^{t}_{min}-y_{t})$ is negative, we demonstrated that their sum is always positive.

Finally, from \eqref{eq:inequality} using \eqref{eq:8}:
\begin{equation*}
\begin{split}
&\dot{\phi}_{k}^{t}-\dot{\phi}_{l}^{t} \leq \bar{w}_{k}- \bar{w}_{l} \\
&-\frac{\delta \cdot K}{n}\Bigl( \sin(\phi^{t}_{max}-y_{t})+\sin(-\phi^{t}_{min}-y_{t})\Bigl) \\
&\leq \bar{w}_{k}- \bar{w}_{l} - \frac{2\delta \cdot K}{n}\sin\Bigl(\frac{D}{2}-\sqrt{\mathcal{E}_{0}-\frac{D^{2}}{2}} \; \Bigl).
\end{split}
\end{equation*}

Thus, $\dot{\phi}_{k}^{t}-\dot{\phi}_{l}^{t} \leq 0$, if  $K \geq \frac{n\cdot |\bar{w}_{k}-\bar{w}_{l}|}{
  2\delta\cdot \sin\bigl( \frac{D}{2}-\sqrt{\mathcal{E}_{0}-\frac{D^{2}}{2}}\bigl)}$.

When $\mathcal{E}_{0}<\frac{D^{2}}{2}$ and condition \eqref{eq:5} is satisfied, $D_{t}$ will be always less than $D$. Indeed, if at time $t$ $D_{t}=D$, then $\phi_{min}^{t}=-\frac{D}{2}+d_{t}$, $\phi_{max}^{t}=\frac{D}{2}+d_{t}$ and $\mathcal{E}_{t} \geq \frac{D^{2}}{2} +2d_{t}^{2} \geq \frac{D^{2}}{2}>\mathcal{E}_{0}$ -- in contradiction to Lemma 2. Therefore, we do not need to have an additional bound \eqref{eq:6} on $K$  to guarantee that $D_{t}\leq D$.
\end{proof}

In Theorem $1$ we have two conditions \eqref{eq:5} and \eqref{eq:6} on the lower bound of the coupling strength $K$, thus the theorem will hold when $K$ satisfies the largest of these two lower bounds.

%


\subsection{Further Refinement Through Optimization}
Similarly to the analytic synchronization condition described in a previous subsection, the optimization-based condition to be introduced in this subsection also guarantees that requirement \eqref{eq:3} of Proposition 1 is satisfied. Numerical techniques, however, allow us to improve the analytic synchronization condition.

There are two bounds on the coupling strength $K$ in Theorem 1, and we will improve bound \eqref{eq:6} using optimization approach. Our optimization approach utilizes additional information that has not been used in the analytic condition, for example, topology information has not been taken into account (except for the minimum nodal degree).
%


For each pair of vertices we solve an optimization problem posed below and find the lower bound on $K$. Then we choose the maximum bound among these obtained $\frac{n(n-1)}{2}$ lower bounds on the coupling strength.


Condition \eqref{eq:3} in Proposition 1 is satisfied if


\begin{equation}\label{eq:9}
K \geq \frac{n\cdot |\bar{w_{k}}-\bar{w_{l}}|}{\sum\limits_{i \in N_{k}}\sin(\phi_{k}-\phi_{i})+\sum \limits_{j \in N_{l}}\sin(\phi_{j}-\phi_{l})}. 
\end{equation}



We find the minimum possible value of the denominator and then obtain corresponding bound on $K$ by \eqref{eq:9}. The  phases of oscillators constitute the phase vector $\vec{\phi}=[\phi_{1}, \dots, \phi_{n}]^{T}$ and are the variables of the optimization problem. The optimization problem is formulated as follows:
\begin{equation*}
\begin{split}
\text{{\bf{minimize}} \hspace{10pt}} &\sum\limits_{i \in N_{k}}\sin(\phi_{k}-\phi_{i})+\sum \limits_{j \in N_{l}}\sin(\phi_{j}-\phi_{l}) \\
\text{{\bf{subject to}} \hspace{7pt}} &\phi_{k}=\phi_{l}+D, \\
\text{ } &\sum\limits_{m=1}^{n}\phi_{m}=0,\\
&\sum\limits_{m=1}^{n}\phi_{m}^{2} \leq\sum\limits_{m=1}^{n}(\phi_{m}^{0})^{2} =\mathcal{E}_{0} ,\\
&\phi_{l}\leq \phi_{m} \leq \phi_{k}, \forall m =1, \dots, n.
\end{split}
\end{equation*}

The first constraint guarantees that the phase distance between oscillators $l$ and $k$ is exactly $D$. Second constraint requires that the sum of all phases is equal to zero. Third condition is necessary because function $\mathcal{E}(t)$ does not exceed its initial value when \eqref{eq:5} is satisfied. The last constraint asks for all the phases to be between $\phi_{l}$ and $\phi_{k}$.

Let $K_{kl}^{*}$ denote the value of the coupling strength in \eqref{eq:9} found with optimization for oscillators $k$ and $l$, and suppose that $K^{*}=\max\limits_{k,l}K_{kl}^{*}$ is the maximum found coupling value among all pairs of oscillators. Then if $K\geq K^{*}$ and condition \eqref{eq:5} is satisfied, system \eqref{eq:2} achieves frequency synchronization. We summarize this result in the following theorem:

{\bf{Theorem 2}} {\it{If $D$ is a constant satisfying $0<D_{0}\leq D<\pi$; $\mathcal{E}_{0}<D^{2}$; if $K\geq K^{*}$, where $K^{*}$ is the bound on the coupling strength obtained with optimization, and if $K$ satisfies \eqref{eq:5}, then there exists a PIS defined by $D_{t}\leq D$ for all $t\geq 0$, and system \eqref{eq:2} achieves frequency synchronization.}}

\begin{proof} Conditions $0<D_{0}\leq D<\pi$ and \eqref{eq:5} are required by Lemma 2 which guarantees that $\mathcal{E}_{t}\leq \mathcal{E}_{0}$ whenever $D_{t}=D$. Condition $K\geq K^{*}$ in turn ensures that requirement \eqref{eq:3} is satisfied and $D_{t}$ does not exceed the value of $D$, which implies that there exists a PIS. Condition $\mathcal{E}_{0}<D^{2}$  guarantees that the optimal value of the optimization problem is strictly positive and there exists a finite positive value of $K_{kl}^{*}$ that satisfies \eqref{eq:9} for each pair of oscillators $k$ and $l$.
\end{proof}



In addition to an improvement of condition \eqref{eq:6}, the numerical method has weaker requirements on the initial phases. Optimization approach can be applied when $0<D_{0}\leq D<\pi$ ($0<D_{0}\leq D\leq \frac{\pi}{2}$ in Theorem 1) and when $\mathcal{E}_{0}<D^{2}$ ($\mathcal{E}_{0}<\frac{3}{4}D^{2}$ in Theorem 1).

\section{NUMERICAL EVALUATION}
\begin{table*}

\caption{Synchronization conditions in our comparative analysis}
\begin{center}
    \begin{tabular}{ | c | c | c |}
    \hline

    & \bf{Bound on Coupling Strength} & \bf{Constraint on Initial Phases} \\ \hline
&&\\
 &$ K \geq \frac{\sigma(\bar{w})\cdot D}{\sqrt{\mathcal{E}_{0}}\cdot L\cdot \sin D}$&\\
    Our condition (optimization approach)  && $\mathcal{E}_{0}< D^{2}<\pi^{2}$ \\
&$K \geq \frac{n\cdot|\bar{w}_{k}-\bar{w}_{l}|}{\sum\limits_{i \in N_{k}}\sin(\phi_{k}-\phi_{i})+\sum \limits_{j \in N_{l}}\sin(\phi_{j}-\phi_{l})}$& \\
 \hline
 && \\
    Condition from \cite{c3}& $K> \frac{2n\cdot\left\|B_{c}^{T}\bar{\bar{w}}\right\|_{2}}{\lambda_{2}\cdot \pi \cdot \text{sinc}(\gamma_{max})  }$ & $\left\|B_{c}^{T}\phi(0)\right\|_{2} < \pi$  \\ 
 && \\
\hline
&&\\
  
   Condition from \cite{c5} &  $K > \frac{\sqrt{2}\sigma(\bar{w})}{L_{*}\cdot\sin D}$  & $\mathcal{E}_{0}<\frac{D^{2}}{2}<\frac{\pi^{2}}{2}$\\

&&\\
    \hline
     \end{tabular}

\label{table_example}
\end{center}
\end{table*}

We compared our results with two existing frequency synchronization conditions. We do not consider here the type of bounds of  \cite{c24}, \cite{c25} and some results of \cite{c3}, since they only provide existence and local stability. To the best of our knowledge the only two bounds that guarantee the existence of a positively invariant set for arbitrary topologies are Theorem 4.6 from \cite{c3}, and \cite{c5}.
%
%
%
%
Therefore, here we compare our condition based on optimization approach, condition from Theorem 4.6 in  \cite{c3} and condition from \cite{c5}.

%
%
%

Each of the three synchronization conditions consists of a bound on the coupling strength and constraints on the initial phases of oscillators.
 In particular, all synchronization conditions require that the difference between any two initial phases is less than $\pi$ (i.e. $D_{0}<\pi$). In addition, each synchronization condition has its own special constraint on the initial phases.

The bounds on the coupling strength and corresponding requirements on the initial phases are summarized in Table 1. Our bound is the maximum between bound \eqref{eq:5} and the optimization bound.
The latter in turn is defined as the maximum bound \eqref{eq:9} among all pairs of oscillators. In our optimization-based synchronization condition, $D$ is a constant whose value can be chosen from the interval  $[D_{0}, \pi)$. While larger values of $D$ make the constraint on the initial phases less restrictive, the bound on the coupling strength is less restrictive on average for smaller values of $D$. In the simulations we assigned a value of $D_{0}$ to the constant $D$. 
In the bound from \cite{c3},  $\lambda_{2}$ is the algebraic connectivity of a given graph, $B_{c} \in \mathbb{R}^{n\times n(n-1)/2}$ is the incidence matrix of the complete graph, $\bar{w}$ is a vector of frequencies, $\phi(0)$ ~--~ vector of initial phases and  $\gamma_{max}=\max\bigl(\frac{\pi}{2},\left\|B_{c}^{T}\phi(0)\right\|_{2}\bigl)$.
 In the bound from \cite{c5}, $D$ is a constant whose value is defined as $\max\bigl(\frac{\pi}{2}, \sqrt{2\mathcal{E}_{0}}\bigl)$, and $L_{*}$ is defined as $L_{*}\triangleq \frac{1}{1+diam(G)\cdot |E^{c}(G)|}$, where $diam(G)$ is the diameter of a graph $G$ and $|E^{c}(G)|$ is the cardinality of its set $E^{c}(G)$.


%

%

In our analysis we compare the requirements on both, the initial phases, and the bounds on the coupling strength.

 \begin{figure}[b]
      \centering
      \includegraphics[scale=0.47]{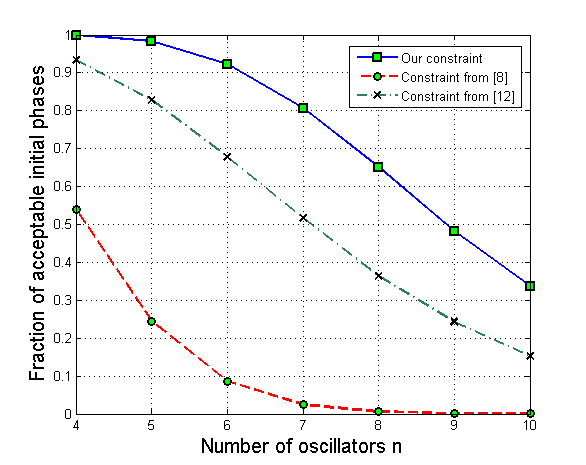}
      \caption{Fractions of random samples of initial phases that satisfy initial phase constraints}
      \label{experiment1}
   \end{figure}

{\bf{Experiment 1 (comparison of the constraints on initial phases).}} 

The first part of our comparison analysis examines how restrictive the constraints on the initial phases are. We created $10^{5}$ samples of initial phases by choosing the phases randomly from the $(0, \pi)$ interval and then subtracting from each sample its mean. The sum of phases in each sample therefore is equal to zero and maximum pairwise phase difference is not greater than $\pi$. We did this for $n=4, \dots, 10$, where $n$ is the number of oscillators in a network or the number of entries in each sample. We then checked for each sample of initial phases, if it satisfies the initial phase constraints. The results are shown on Fig. \ref{experiment1}, where the x-axis corresponds to the number of oscillators $n$ and the y-axis corresponds to a fraction of samples that satisfy the initial phase constraints for each synchronization condition. From this figure we can see that our constraint on the initial phases is the least restrictive, whereas the constraint corresponding to condition from \cite{c3} is the most restrictive on average. From this result we may conclude that our synchronization condition has a larger region of applicability when initial phases are randomly chosen from the $(0, \pi)$ interval.

{\bf{Experiment 2 (comparison of the bounds on coupling strength).}}

 \begin{figure*}[thpb]
\label{experiment2}
\begin{center}
      \centering
\begin{subfigure}[t]{0.3\textwidth}
      \includegraphics[width=\textwidth]{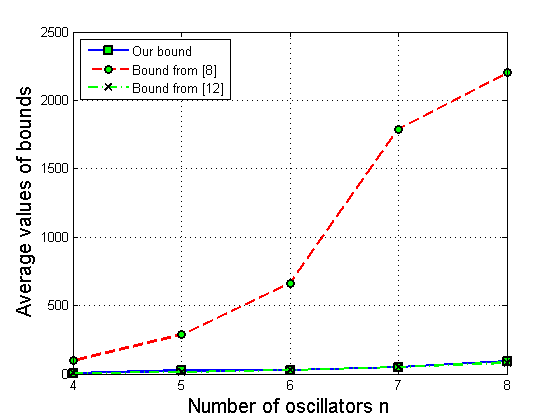}
      \caption{Average bounds, chain topology }      
\end{subfigure}
\begin{subfigure}[t]{0.3\textwidth}
      \includegraphics[width=\textwidth]{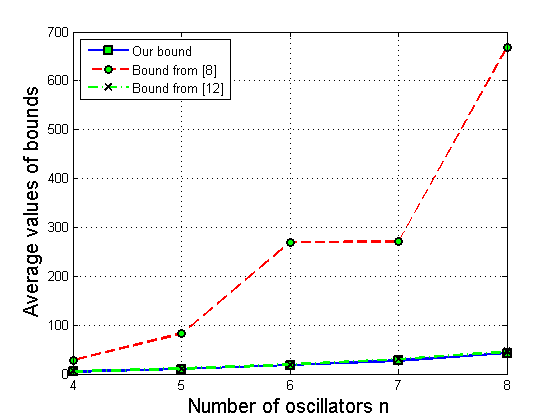}
      \caption{Average bounds, star-tree topology}      
\label{4b}
\end{subfigure}
\begin{subfigure}[t]{0.3\textwidth}
      \includegraphics[width=\textwidth]{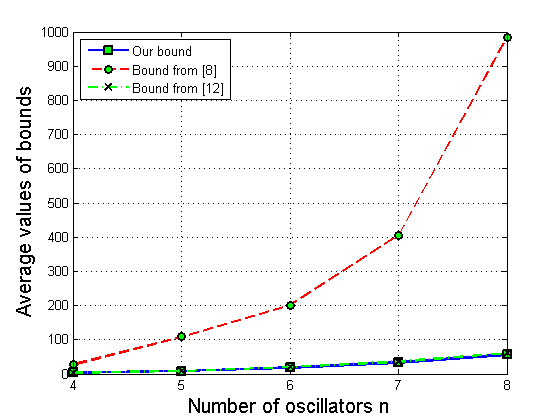}
      \caption{Average bounds, ring topology }     
\end{subfigure}

\begin{subfigure}[t]{0.3\textwidth}
      \includegraphics[width=\textwidth]{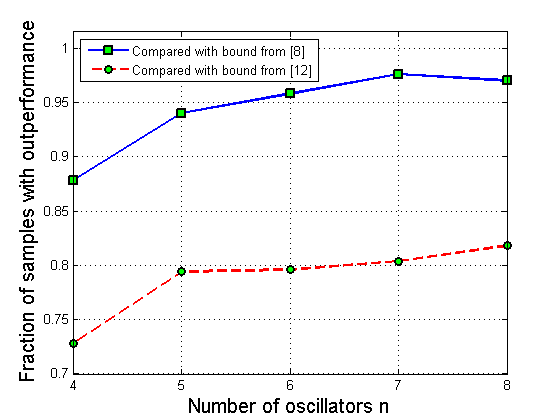}
      \caption{Fractions of samples with outperformance, chain topology }     
\end{subfigure}
\begin{subfigure}[t]{0.3\textwidth}
      \includegraphics[width=\textwidth]{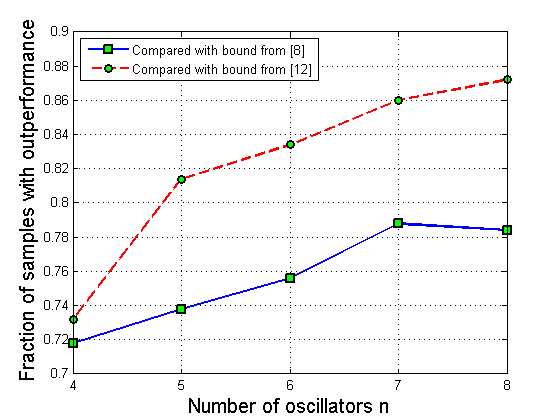}
      \caption{Fractions of samples with outperformance, star-tree topology}     
\end{subfigure}
\begin{subfigure}[t]{0.3\textwidth}
      \includegraphics[width=\textwidth]{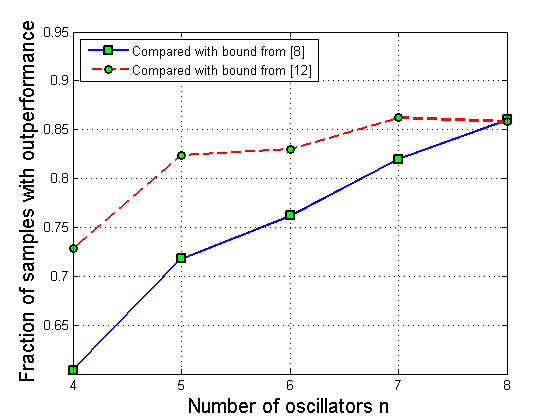}
      \caption{Fractions of samples with outperformance, ring topology }     
\end{subfigure}
 \label{experiment2}
\caption{Average values of bounds (a-c) and fractions of samples on which our bound outperforms other bounds (from \cite{c3} and \cite{c5}), (d-f) for chain, star-tree and ring topologies.}
\end{center}
 \label{experiment2}
   \end{figure*}

In the second part of our comparison we applied the bounds on the coupling strength only to the samples of the initial phases that satisfy all three initial phase constraints. We performed the comparison for star-tree (a tree with one node connected to all others), ring and chain topologies with the number of oscillators $n$ varying from 4 to 8. Because all three bounds contain oscillators' frequencies, we also had to create random samples of frequencies. The simulation process was organized as follows: we first fixed the underlying topology, that is the topology type (star-tree, ring or chain) together with the number of oscillators $n$. After that, for each topology we generated $500$ pairs of random samples of frequencies and initial phases.
Then, for each pair of samples of  frequencies and initial phases we calculated values of three bounds on the coupling strength.
We then found average values of the three bounds over $500$ samples for each topology. We also calculated for each topology the fraction of samples for which our bound outperforms bounds from \cite{c3} and from \cite{c5}.

While all three bounds can be applied only when the initial maximum phase difference $D_{0}$ is from the interval $(0,\pi)$, there is no similar requirement on the frequencies. To figure out if an interval from which the frequencies are sampled influences relative performance of the three bounds, we performed experiment 2 for two sample intervals of the frequencies: $(0,1)$ and $(0,10)$. Because results for these two cases are very similar, we provide only the results corresponding to the interval $(0,1)$.

The average values of bounds are plotted on Fig. 3 (a-c). 
The performances of our bound and bound from \cite{c5} seem to be very similar, especially when plotted together with the performance of bound from \cite{c3}. However, on Fig. 3 (d-f)  we plotted for each topology the fractions of frequency-phases samples for which bounds from \cite{c3} and from \cite{c5} are more restrictive. As we can see, on the majority of samples our bound outperforms the two other bounds.

%
%

We can now summarize the results of our comparison analysis of the three synchronization conditions. The constraint on the initial phases for our condition seems to be the least restrictive  compared to the similar constraints for bounds two other bounds. In addition, on the majority of samples our coupling strength bound based on optimization approach is the least restrictive as can be seen on Fig. 3.

%
%

We used Matlab's R2012a $\it{GlobalSearch}$ function from the Global Optimization Toolbox with the following options: MaxFunEvals = 300000, MaxIter = 500000, TolFun = $10^{-10}$, TolCon = $10^{-10}$, TolX = $10^{-10}$.

\section{CONCLUSION}

This paper studies synchronization of heterogeneous Kuramoto oscillators with arbitrary underlying topology. We provide novel sufficient conditions on the coupling strength that guarantee the existence of a Positively Invariant Set (PIS) and then use LaSalle's Invariance Principle to show frequency synchronization. Moreover, we provide an optimization framework that can further improve our bounds. We illustrate these results with simulations performed for chain, ring and star-tree topologies. Our bounds consistently improve existing bounds on average for every investigated case.

\section*{APPENDIX}

{\bf{Lemma 1}} {\it{If $\sum\limits_{i=1}^{n}\phi_{i}^{t}=0$, then

 $$ L \cdot n \cdot \mathcal{E}_{t} \leq \sum\limits_{(i,j) \in E}|\phi_{i}^{t}-\phi_{j}^{t}|^{2} \leq n\cdot \mathcal{E}_{t}, $$
where  $L \triangleq \frac{1}{1+\sum\limits_{(k,l) \in E^{c}}dist(k,l)}$.}}

\begin{proof}
For simplicity, in the proof we will omit the time transcript $t$.
Upper bound:
\begin{equation*}
\begin{split}
\sum\limits_{(i,j) \in E}|\phi_{i}-\phi_{j}|^{2} &\leq \frac{1}{2}\sum\limits_{k,l =1}^{n}|\phi_{k}-\phi_{l}|^{2}\\
&=\frac{1}{2}\sum\limits_{k,l=1}^{n}(\phi_{k}^{2}+\phi_{l}^{2}-2\phi_{k}\phi_{l}) =n\cdot \mathcal{E}.
\end{split}
\end{equation*}

Lower bound: if edge $(k \rightarrow l) \in E^{c}$, then by a triangle inequality:
$$|\phi_{k}-\phi_{l}|^{2} \leq dist(k,l)\cdot \sum\limits_{(i,j) \in E}|\phi_{i}-\phi_{j}|^{2}.$$
Thus:
$$\sum\limits_{(k,l) \in E^{c}}|\phi_{k}-\phi_{l}|^{2}\leq \Bigl(\sum\limits_{(k,l)\in E^{c}}dist(k,l)\Bigl)\cdot \Bigl(\sum\limits_{(i,j) \in E}|\phi_{i}-\phi_{j}|^{2}\Bigl).$$
Then
\begin{equation*}
\begin{split}
\sum\limits_{k,l =1}^{n}|\phi_{k}&-\phi_{l}|^{2}=2\sum\limits_{(k,l) \in E}|\phi_{k}-\phi_{l}|^{2} + 2\sum\limits_{(k,l) \in E^{c}}|\phi_{k}-\phi_{l}|^{2} \\
&\leq 2\Bigl(1+\sum\limits_{(k,l)\in E^{c}}dist(k,l)\Bigl)\cdot \sum\limits_{(i,j) \in E}|\phi_{i}-\phi_{j}|^{2}.
\end{split}
\end{equation*}
Therefore:
$$\sum\limits_{(i,j) \in E}|\phi_{i}-\phi_{j}|^{2}\geq \frac{L}{2}\cdot \sum\limits_{k,l=1 }^{n}|\phi_{k}-\phi_{l}|^{2}=L\cdot n\cdot \mathcal{E}.$$
\end{proof}


{\bf{Lemma 2}} {\it{ If $D_{t} \leq D<\pi \; \forall t \in [0, T]$, then function $\mathcal{E}$ satisfies the following differential inequality on $[0,T]$:
$$
\frac{d}{dt}\mathcal{E} \leq 2\sigma(\bar{w})\cdot  \sqrt{\mathcal{E}}-2K\cdot L\cdot \Bigl(\frac{\sin D}{D}\Bigl)\cdot \mathcal{E}, 
$$
and if in addition
$$ K\geq \frac{\sigma(\bar{w})\cdot D}{\sqrt{\mathcal{E}_{0}}\cdot L\cdot \sin D} \;,$$
then $\mathcal{E}_{t}$ will be upper bounded by $\mathcal{E}_{0}$: $\mathcal{E}_{t} \leq \mathcal{E}_{0} \; \forall t\in [0, T]$.}}
\begin{proof}
Multiplying $i^{th}$ equation of \eqref{eq:2} by $2\phi_{i} \; (\forall i=1,\dots ,n)$ and summing them together:
$$\frac{d}{dt}\sum\limits_{i=1}^{n}\phi_{i}^{2}=2\sum\limits_{i=1}^{n}\bar{w}_{i}\phi_{i}-\frac{K}{n}\sum\limits_{i,j=1}^{n}(\phi_{j}-\phi_{i})\sin(\phi_{j}-\phi_{i}).$$
It can be verified that:
$$\frac{\sin x}{x} \geq \frac{\sin D}{D}, \; \forall x\in [-D, D].$$
Thus, using Lemma 1 and the Schwartz inequality:
\begin{equation*}
\begin{split}
\frac{d}{dt}\mathcal{E} &\leq 2\sigma(\bar{w})\cdot \sqrt{\mathcal{E}}-\frac{2K}{n}\cdot \Bigl(\frac{\sin D}{D}\Bigl)\cdot \sum\limits_{(i,j)\in E}(\phi_{i}-\phi_{j})^{2} \\
&\leq 2\sigma(\bar{w})\cdot \sqrt{\mathcal{E}}-2K\cdot L\cdot \Bigl(\frac{\sin D}{D}\Bigl)\cdot  \mathcal{E}.
\end{split}
\end{equation*}
We now consider the following differential equation:
\begin{equation*}
\begin{split}
&\frac{dz}{dt}=\sigma(\bar{w})-K\cdot L\cdot \Bigl(\frac{\sin D}{D}\Bigl) \cdot z,\\
&z(0)=\sqrt{\mathcal{E}_{0}}.
\end{split}
\end{equation*}
This equation has one asymptotically stable equilibrium $z_{e}$, and $z(t)$ monotonically decreases to $z_{e}$ if condition \eqref{eq:5} is satisfied and therefore $z(t)\leq z(0) \; \forall t\geq 0$. By comparison principle, 
$$\sqrt{\mathcal{E}(t)}\leq z(t), t \in [0, T], $$
and thus $\mathcal{E}_{t} \leq \mathcal{E}_{0} \; \forall t \in [0,T]$.
\end{proof}


\end{document}